\pdfoutput=1
\RequirePackage{ifpdf}
\ifpdf
\documentclass[pdftex]{sigma}
\else
\documentclass{sigma}
\fi

\numberwithin{equation}{section}

\numberwithin{theorem}{section}
\numberwithin{proposition}{section}
\numberwithin{lemma}{section}
\numberwithin{corollary}{section}
\numberwithin{definition}{section}
\numberwithin{example}{section}
\numberwithin{remark}{section}
\numberwithin{note}{section}

\usepackage[all]{xy}

\theoremstyle{remark}

\newcommand{\mR}{\mathbb{R}}
\newcommand{\mC}{\mathbb{C}}
\newcommand{\mN}{\mathbb{N}}
\newcommand{\mE}{\mathbb{E}}
\newcommand{\mZ}{\mathbb{Z}}
\newcommand{\mS}{\mathbb{S}}
\newcommand{\mV}{\mathbb{V}}

\newcommand{\cD}{\mathcal{D}}
\newcommand{\cM}{\mathcal{M}}
\newcommand{\cH}{\mathcal{H}}

\newcommand{\cF}{\mathcal{F}}
\newcommand{\cP}{\mathcal{P}}
\newcommand{\cC}{\mathcal{C}}

\newcommand{\cS}{\mathcal{S}}
\newcommand{\cG}{\mathcal{G}}
\newcommand{\cL}{\mathcal{L}}

\newcommand{\dD}{\textbf{D}}

\newcommand{\ux}{\underline{x}}
\newcommand{\uz}{\underline{z}}

\newcommand{\uy}{\underline{y}}

\newcommand{\upx}{\partial_{\underline{x}}}

\def\a{{\alpha}}
\def\b{{\beta}}
\def\g{{\gamma}}
\def\k{{\kappa}}

\def\l{{\lambda}}
\def\d{{\delta}}
\def\o{{\omega}}
\def\s{{\sigma}}
\def\la{{\langle}}
\def\ra{{\rangle}}

\begin{document}
\allowdisplaybreaks

\renewcommand{\PaperNumber}{010}

\FirstPageHeading

\ShortArticleName{The Clif\/ford Deformation of the Hermite Semigroup}

\ArticleName{The Clif\/ford Deformation of the Hermite Semigroup}

\Author{Hendrik DE BIE~$^\dag$, Bent {\O}RSTED~$^\ddag$, Petr SOMBERG~$^\S$ and Vladimir
SOU{\v{C}}EK~$^\S$}

\AuthorNameForHeading{H.~De Bie, B.~{\O}rsted, P.~Somberg and V.~Sou{\v{c}}ek}

\Address{$^\dag$~Department of Mathematical Analysis, Ghent University, Galglaan 2, 9000 Gent,
Belgium}
\EmailD{\href{mailto:Hendrik.DeBie@UGent.be}{Hendrik.DeBie@UGent.be}}

\Address{$^\ddag$~Department of Mathematical Sciences, University of Aarhus,\\
\hphantom{$^\ddag$}~Building 530, Ny Munkegade, DK 8000, Aarhus C, Denmark}
\EmailD{\href{mailto:orsted@imf.au.dk}{orsted@imf.au.dk}}

\Address{$^\S$ Mathematical Institute of Charles University, Sokolovsk\'a 83, 186 75 Praha, Czech
Republic}
\EmailD{\href{mailto:somberg@karlin.mff.cuni.cz}{somberg@karlin.mff.cuni.cz},
\href{mailto:soucek@karlin.mff.cuni.cz}{soucek@karlin.mff.cuni.cz}}

\ArticleDates{Received September 21, 2012, in f\/inal form January 29, 2013; Published online February 05, 2013}

\Abstract{This paper is a~continuation of the paper [De~Bie H., {\O}rsted B., Somberg P.,
Sou{\v{c}}ek V., \textit{Trans.
Amer.
Math.
Soc.} \textbf{364} (2012), 3875--3902], investigating a~natural radial
deformation of the Fourier transform in the setting of Clif\/ford analysis.
At the same time, it gives
extensions of many results obtained in [Ben~Sa{\"{\i}}d S., Kobayashi T., {\O}rsted B.,
\textit{Compos.
Math.} \textbf{148} (2012), 1265--1336].
We establish the analogues of Bochner's formula and the Heisenberg uncertainty relation in the
framework of the (holomorphic) Hermite semigroup, and also give a~detailed analytic treatment of
the series expansion of the associated integral transform.
}


\Classification{33C52; 30G35; 43A32}

\section{Introduction}
\label{intro}

It is well-known that the classical Dirac operator and its Fourier symbol generate via Clif\/ford
multiplication a~natural Lie superalgebra $\mathfrak{osp}(1|2)$ contained in the Clif\/ford--Weyl
algebra.
More surprisingly, this carries over to a~natural family of deformations
of the Dirac operator, see~\cite{H12}.
Moreover, it is possible to def\/ine a~Fourier transform naturally associated to the deformed family.

The novelty of the present article is that we let group theory be the
guiding principle in def\/ining operators and transformations,
in the next step followed by a~study of explicit (analytic) properties
for naturally arising eigenfunctions and kernel functions.
Thus the main aim is to f\/ind the kernel function for the Fourier transform connected with our
deformation, and also to study its associated holomorphic semigroup
regarded as a~particular descendant of the Gelfand--Gindikin program analyzing
representations of reductive Lie groups, see, e.g.,~\cite{Ol} and the discussion in~\cite{Orsted2}.

Let us now recall the basic setup and results from~\cite{H12}
and also discuss further aspects of our construction.
The deformation family of Dunkl--Dirac operators
\begin{gather*}
\dD=r^{1-\frac{a}{2}}\cD_{\k}+b r^{-\frac{a}{2}-1}\ux+c r^{-\frac{a}{2}-1}\ux\mE,\qquad a,b,c\in\mR,
\end{gather*}
together with the radial deformation of the coordinate function
\begin{gather*}
\ux_{a}=r^{\frac{a}{2}-1}\ux,\qquad r=\sqrt{\sum_{i=1}^{m}x_{i}^{2}},
\end{gather*}
forms a~realization of $\mathfrak{osp}(1|2)$ in the Clif\/ford--Weyl algebra.
Here $\cD_{\k}=\sum\limits_{i=1}^{m}e_{i}T_{i}$ with $T_{i}$ the Dunkl operators, $\ux=
\sum\limits_{i=1}^{m}e_{i}x_{i}$ and $\mE=\sum\limits_{i=1}^{m}x_{i}\partial_{x_{i}}$.
The $e_{i}$ are generators of the Clif\/ford algebra~$\cC l_{m}$.
See also the next section for more details.

We will show in Proposition~\ref{equivarianceG} that this realization builds a~Howe dual pair with
$\widetilde{\cG}$. Here the group $\widetilde{\cG}$ is the double cover (contained in the Pin
group) of the f\/inite ref\/lection group $\cG$ used in the construction of the Dunkl operators.

The Fourier transform is then def\/ined by
\begin{gather*}
\cF_{\dD}=e^{i\frac{\pi}{2}\left(\frac{1}{2}+\frac{\mu-1}{a(1+c)}\right)}e^{\frac{-i\pi}{
2a(1+c)^{2}}\left(\dD^{2}-(1+c)^{2}\ux_{a}^{2}\right)},
\end{gather*}
where $L=\dD^{2}-(1+c)^{2}\ux_{a}^{2}$ is the generalized Hamiltonian and $\mu$ the Dunkl
dimension.
The main aim of the present paper is to f\/ind an integral expression for this Fourier transform,
\begin{gather*}
\cF_{\dD}(f)(y)=\int_{\mR^{m}}K(x,y)f(x)h(r_{x})dx
\end{gather*}
with $h(r_{x})dx$ the measure associated to $\dD$ and $K(x,y)$ the integral kernel to be
determined.
Note that this ties in with recent work on generalized Fourier transforms in dif\/ferent contexts,
e.g., analysis on minimal representations of reductive groups (see~\cite{MR2134314, MR2401813,
KM3}) or integral transforms in Clif\/ford analysis (see~\cite{DBNS, DBXu}).

The deformation of the classical Hamiltonian for the harmonic oscillator is visualized in the
following f\/igure:
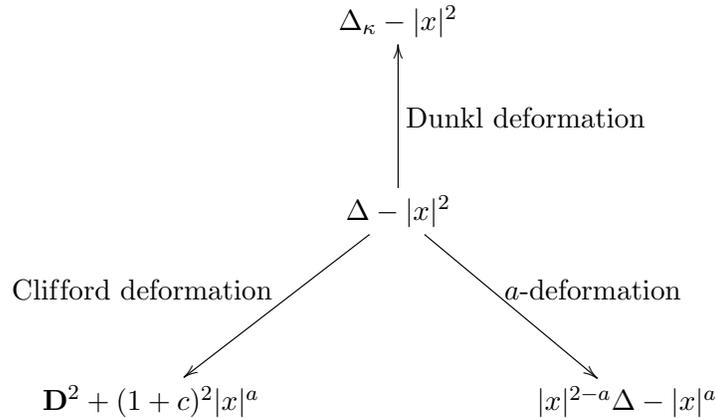
$$
\xymatrix{
&\Delta_{\k}-|x|^{2}&\\
\\
&\Delta-|x|^{2}\ar[uu]_{\mbox{Dunkl deformation}}\ar[ddl]_{\mbox{Clif\/ford
deformation\hphantom{aa}}}\ar[ddr]^{\mbox{$a$-deformation}}&\\
\\
\dD^{2}+(1+c)^{2}|x|^{a}&&|x|^{2-a}\Delta-|x|^{a}
}
$$
The Dunkl deformation is by now quite standard and described for example in~\cite{Du4}.
The $a$-deformation is the subject of the paper~\cite{Orsted2} and is a~scalar radial deformation
of the harmonic oscillator.
Our Clif\/ford deformation is also a~radial deformation but richer in the sense that Clif\/ford
algebra- (or spinor)-valued functions are involved.

In this paper we will thus f\/ind a~series representation of the kernel function for our new Fourier
transform $\cF_{\dD}$, and also study the holomorphic semigroup with generator $L$.
The main results are Theorem~\ref{PropSemigroup} on the operator properties of the semigroup,
Theorem~\ref{StatementsFourier} on the Fourier transform intertwining the Dirac operator and the
Clif\/ford multiplication, Proposition~\ref{Bochner} on the Bochner identities, and
Proposition~\ref{Heisenberg} on the Heisenberg uncertainty relation.
Finally in Theorem~\ref{MasterFormula} we give the analogue of what is sometimes called the
``Master formula'' in the context of Dunkl operators (see, e.g.,~\cite[Lemma~4.5(1)]{MR1620515}
or~\cite{C}).

The paper is organized as follows.
In Section~\ref{section2} we repeat basic notions on Clif\/ford algebras and Dunkl operators needed
in the rest of the paper.
In Section~\ref{section3} we construct intertwining operators to reduce our radially deformed Dirac
operator to its simplest form.
Subsequently, in Section~\ref{section4} we discuss the representation theoretic content of our
deformation and solve the spectral problem of the associated Hamiltonian.
In Section~\ref{section5}, we obtain the reproducing kernels for spaces of spherical monogenics,
which allows us to construct the kernel of the holomorphic semigroup in Section~\ref{section6}.
Section~\ref{section7} contains the results on the (deformed) Fourier transform.
Further properties are collected in Section~\ref{section8}.
Finally, we summarize some results on special functions used in the paper in
Appendix~\ref{appendixA} and give a~list of notations in Appendix~\ref{appendixB}.

\section{Preliminaries}\label{section2}

In this section we collect some basic results on Clif\/ford algebras and Dunkl operators.

\subsection{Clif\/ford algebras}\label{section2.1}
Let $\mV$ be a~vector space of dimension $m$ with a~given negative def\/inite quadratic form and
let~$\cC l_{m}$ be the corresponding Clif\/ford algebra.
If $\{e_i\}$ is an orthonormal basis of $\mV$, then $\cC l_{m}$ is generated by $e_{i}$, $i=1,
\ldots,m$, with the relations
\begin{gather*} 
e_{i} e_{j} + e_{i} e_{j} = 0, \quad i \neq j,\qquad
 e_{i}^{2} = -1.
\end{gather*}
The algebra $\cC l_{m}$ has dimension $2^{m}$ as a~vector space over $\mR$.
It can be decomposed as $\cC l_{m}=\oplus_{k=0}^{m}\cC l_{m}^{k}$
with $\cC l_{m}^{k}$ the space of $k$-vectors def\/ined by
\begin{gather*}
\cC l_{m}^{k}:=\text{span}\{e_{i_{1}}\cdots e_{i_{k}},\; i_{1}<\cdots<i_{k}\}.
\end{gather*}
The projection on the space of $k$-vectors is denoted by $[\cdot]_{k}$.

The operator $\bar{.}$ is the main anti-involution on the Clif\/ford algebra $\cC l_{m}$ def\/ined by
\begin{gather*}
\overline{a b}=\overline{b}\overline{a},\qquad\overline{e_{i}}=-e_{i},\quad i=1,\ldots,m.
\end{gather*}
Similarly we have the automorphism $\epsilon$ given by
\begin{gather*}
\epsilon(a b)=\epsilon(a)\epsilon(b),\qquad\epsilon(e_{i})=-e_{i},\quad i=1,\ldots,m.
\end{gather*}

In the sequel, we will always consider functions $f$ taking values in $\cC l_{m}$, unless
explicitly mentioned.
Such functions can be decomposed as
\begin{gather*}
f(x)=f_{0}(x)+\sum_{i=1}^{m}e_{i}f_{i}(x)+\sum_{i<j}e_{i}e_{j}f_{ij}(x)+\cdots+e_{1}
\cdots e_{m}f_{1\ldots m}(x)
\end{gather*}
with $f_{0},f_{i},f_{ij},\ldots,f_{1\ldots m}$ all real-valued functions.

Several important groups can be embedded in the Clif\/ford algebra.
Note that the space of $1$-vectors in $\cC l_{m}$ is canonically isomorphic to $\mV\cong\mR^{m}$.
Hence we can def\/ine
\begin{gather*}
\mathrm{Pin}(m)=\left\{s_{1}s_{2}\cdots s_{n}\,|\,n\in\mN,s_{i}\in\cC l_{m}^{1}\ \text{such that}\ s_{i}^{2}=-1\right\},
\end{gather*}
i.e., the Pin group is the group of products of unit vectors in~$\cC l_{m}$.
This group is a~double cover of the orthogonal group~${\rm O}(m)$ with covering map $p:\mathrm{Pin}(m)
\rightarrow {\rm O}(m)$, which we will describe explicitly in the next section.

Similarly we def\/ine
\begin{gather*}
\mathrm{Spin}(m)=\left\{s_{1}s_{2}\cdots s_{2n}\,|\,n\in\mN,s_{i}\in\cC l_{m}^{1} \ \text{such that} \ s_{i}^{2}=-1\right\},
\end{gather*}
i.e., the Spin group is the group of even products of unit vectors in $\cC l_{m}$.
This group is a~double cover of ${\rm SO}(m)$.
For more information about Clif\/ford algebras and analysis, we refer the reader to~\cite{MR1169463,MR1130821}.

\subsection{Dunkl operators}\label{section2.2}
Denote by $\langle\cdot,\cdot
\rangle$ the standard Euclidean scalar product in $\mR^{m}$ and by $|x|=\langle x,
x\rangle^{1/2}$ the associated norm.
For $\alpha\in\mR^{m}\setminus \{0\}$, the ref\/lection $r_{\alpha}$ in the hyperplane orthogonal to
$\alpha$ is given by
\begin{gather*}
r_{\alpha}(x)=x-2\frac{\langle\alpha,x\rangle}{|\alpha|^{2}}\alpha,\qquad x\in\mR^{m}.
\end{gather*}

A root system is a~f\/inite subset $R\subset\mR^{m}$ of non-zero vectors such that, for every
$\alpha\in R$, the associated ref\/lection $r_{\alpha}$ preserves $R$.
We will assume that $R$ is reduced, i.e.
$R\cap\mR\alpha=\{\pm\alpha\}$ for all $\alpha\in R$.
Each root system can be written as a~disjoint union $R=R_{+}\cup(-R_{+})$, where $R_{+}$ and~$-R_{+}$ are separated by a~hyperplane through the origin.
$R_+$ is called a~positive subsystem of the root system $R$.
The subgroup $\cG\subset {\rm O}(m)$ generated by the ref\/lections $\{r_{\alpha}|\alpha\in R\}$ is
called the f\/inite ref\/lection group associated with~$R$.
We will also assume that $R$ is normalized such that $\langle\alpha,\alpha\rangle=2$ for all
$\alpha\in R$.
For more information on f\/inite ref\/lection groups we refer the reader to~\cite{Humph}.

If we identify $\alpha$ with a~$1$-vector in $\cC l_{m}$ (and hence $\alpha/\sqrt{2}$ with an
element in $\mathrm{Pin}(m)$), we can rewrite the ref\/lection $r_{\alpha}$ as
\begin{gather*}
r_{\alpha}(x)=\frac{1}{2}\alpha\ux\alpha
\end{gather*}
with $\ux=\sum\limits_{i=1}^{m}e_{i}x_{i}$.
Generalizing this map gives us the covering map $p$ from $\mathrm{Pin}(m)$ to ${\rm O}(m)$ as
\begin{gather*}
p(s)(x)=\epsilon(s)\ux s^{-1},\qquad s\in \mathrm{Pin}(m).
\end{gather*}
In particular, we obtain a~double cover of the ref\/lection group $\cG$ as $\widetilde{\cG}=
p^{-1}(\cG)$ (see also the discussion in~\cite{BCT}).

A multiplicity function $\kappa$ on the root system $R$ is a~$\cG$-invariant function $\kappa:R
\rightarrow\mC$, i.e.
$\k(\alpha)=\k(h\alpha)$ for all $h\in\cG$.
We will denote $\k(\alpha)$ by $\k_{\alpha}$.
We will always assume that the multiplicity function is real and satisf\/ies $\k\geq0$.
This assumption is, e.g., 
necessary to obtain the subsequent formula~\eqref{skewnessDunkl}, which is crucial for the sequel.

Fixing a~positive subsystem $R_{+}$ of the root system~$R$ and a~multiplicity function~$\k$, we
introduce the Dunkl operators~$T_{i}$ associated to~$R_{+}$ and $\k$ by (see~\cite{MR951883,
MR1827871})
\begin{gather*}
T_{i}f(x)=\partial_{x_{i}}f(x)+\sum_{\alpha\in R_{+}}\k_{\alpha}\alpha_{i}\frac{f(x)-
f(r_{\alpha}(x))}{\langle\alpha,x\rangle},\qquad f\in C^{1}(\mR^{m}).
\end{gather*}
An important property of the Dunkl operators is that they commute, i.e.
$T_{i}T_{j}=T_{j}T_{i}$.

The Dunkl Laplacian is given by $\Delta_{\k}=\sum\limits_{i=1}^{m}T_i^2$, or more explicitly by
\begin{gather*}
\Delta_{\k}f(x)=\Delta f(x)+2\sum_{\alpha\in R_{+}}\k_{\alpha}\left(\frac{\langle\nabla
f(x),\alpha\rangle}{\langle\alpha,x\rangle}-\frac{f(x)-f(r_{\alpha}(x))}{\langle\alpha,
x\rangle^{2}}\right)
\end{gather*}
with $\Delta$ the classical Laplacian and $\nabla$ the gradient operator.
We also def\/ine the constant
\begin{gather*}
\mu=\frac{1}{2}\Delta_{\k}|x|^2=m+2\sum_{\alpha\in R_+}\k_{\alpha},
\end{gather*}
called the Dunkl-dimension.

It is possible to construct an intertwining operator $V_{\k}$ connecting the classical derivatives
$\partial_{x_{j}}$ with the Dunkl operators $T_{j}$ such that $T_{j}V_{\k}=V_{\k}
\partial_{x_{j}}$ (see, e.g.,~\cite{MR1273532}).
Note that explicit formulae for~$V_{\k}$ are only known in a~few special cases.

The weight function related to the root system $R$ and the multiplicity function $\k$ is given by
$w_{\k}(x)=\prod\limits_{\alpha\in R_{+}}|\langle\alpha,x\rangle|^{2\k_{\alpha}}$.
For suitably chosen functions $f$ and $g$ one then has the following property of integration by
parts (see~\cite{Du4})
\begin{gather}
\int_{\mR^{m}}(T_{i}f)g w_{\k}(x)dx=-\int_{\mR^{m}}f\left(T_{i}g\right)w_{\k}(x)dx.
\label{skewnessDunkl}
\end{gather}
For more information about Dunkl operators we refer the reader to~\cite{MR1827871, MR2022853}.

The starting point in the subsequent analysis is the Dunkl--Dirac operator, given by
\begin{gather*}
\cD_{\k}=\sum_{i=1}^{m}e_{i}T_{i}.
\end{gather*}
Together with the vector variable $\ux=\sum\limits_{i=1}^{m}e_{i}x_{i}$ this Dunkl--Dirac operator
generates a~copy of~$\mathfrak{osp}(1|2)$, see~\cite{Orsted} or the subsequent Theorem~\ref{ospFamily}.
In particular, we have
\begin{gather*}
\cD_{\k}^{2}=-\Delta_{\k}\qquad\text{and}\qquad\ux^{2}=-|\ux|^{2}=-r^{2}=-
\sum_{i=1}^{m}x_{i}^{2}.
\end{gather*}

\section{Intertwining operators}\label{section3}

Let, for $a,b\in\mR$, $P$ and $Q$ be two operators def\/ined by
\begin{gather*}
P f(\ux)=r^{b}f\left(\left(\frac{a}{2}\right)^{\frac{1}{a}}\ux r^{\frac{2}{a}-1}\right),\qquad
Q f(\ux)=r^{-\frac{ab}{2}}f\left(\left(\frac{2}{a}\right)^{\frac{1}{2}}\ux r^{\frac{a}{2}-1}\right).
\end{gather*}
These two operators act as generalized Kelvin transformations.
Indeed, one can easily compute their composition
\begin{gather*}
Q P=P Q=\left(\frac{2}{a}\right)^{\frac{b}{2}}.
\end{gather*}
We will show that these operators allow to reduce the Dirac operator $\dD$ to a~simpler form.

We have the following proposition, where $\mE=\sum\limits_{i=1}^{m}x_{i}\partial_{x_{i}}$ denotes the
Euler operator.
Recall also from the introduction that $\ux_{a}=r^{\frac{a}{2}-1}\ux$.

\begin{proposition}
One has the following intertwining relations
\begin{gather*}
\left(\frac{a}{2} \right)^{\frac{b-1}{2}} Q \left( \cD_{\k} + b r^{-2} \ux + c r^{-2}\ux \mE
\right) P = r^{1- \frac{a}{2}}\cD_{k} + \b r^{-\frac{a}{2}-1} \ux + \g r^{-\frac{a}{2}-1}\ux \mE,\\
\left(\frac{a}{2} \right)^{\frac{b+1}{2}} Q \ux P = \ux_{a}
\end{gather*}
with
$
\b= 2b+ bc$, $
\g= \frac{2}{a} (1+c) -1$.
\end{proposition}

\begin{proof}
In~\cite[Proposition 3]{H12}, we already proved that
\begin{gather*}
\left(\frac{a}{2}\right)^{\frac{b-1}{2}}Q\left(\cD_{\k}\right)P=r^{1-\frac{a}{2}}\cD_{k}+
b r^{-\frac{a}{2}-1}\ux+\left(\frac{2}{a}-1\right)r^{-\frac{a}{2}-1}\ux\mE=\ux_{a}.
\end{gather*}
Similarly we obtain
\begin{gather*}
\left(\frac{a}{2}\right)^{\frac{b-1}{2}}Q\left(r^{-2}\ux\right)P=r^{-\frac{a}{2}-1}\ux
\end{gather*}
and
\begin{gather*}
\left(\frac{a}{2}\right)^{\frac{b-1}{2}}Q\left(r^{-2}\ux\mE\right)P=b r^{-\frac{a}{2}-1}
\ux+\left(\frac{2}{a}\right)r^{-\frac{a}{2}-1}\ux\mE.
\end{gather*}
This completes the proof of the proposition.
\end{proof}

So we are reduced to the study of the operator
\begin{gather*}
\dD=\cD_{\k}+b r^{-2}\ux+c r^{-2}\ux\mE,
\end{gather*}
where $b,c\in\mR$, $c\neq-1$.
Here, the term $b r^{-2}\ux$ can also be removed.
Indeed, we have
\begin{gather*}
r^{-\alpha}\left(\cD_{\k}+b r^{-2}\ux+c r^{-2}\ux\mE\right)r^{\alpha}=\cD_{\k}+c
r^{-2}\ux\mE,
\end{gather*}
when $\alpha=-b/(1+c)$.

As a~result of the previous discussion, we see that it is suf\/f\/icient to study the function theory
for the operator
\begin{gather*}
\dD=\cD_{\k}+c r^{-2}\ux\mE,
\end{gather*}
where we have put $a=2$, $b=0$.
Furthermore, we will restrict ourselves to the case $c>-1$ for reasons that will become clear in
Proposition~\ref{partIntProp}.
Similarly, we no longer need to consider $\ux_{a}$ but can restrict ourselves to $\ux$.
Now we repeat the basic facts concerning this operator we need in the sequel.
All the results are taken from~\cite{H12}, putting $a=2$, $b=0$.

\begin{theorem}\label{ospFamily}
The operators $\dD$ and $\ux$ generate a~Lie superalgebra, isomorphic to $\mathfrak{osp}(1|2)$,
with the following relations
\begin{alignat*}{3}
& \{\ux,\dD\}=-2(1+c)\left(\mE+\frac{\delta}{2}\right),\qquad&&
\left[\mE+\frac{\delta}{2},\dD\right]=-\dD,& \\
& \left[\ux^{2},\dD\right]=2(1+c)\ux,\qquad&&
\left[\mE+\frac{\delta}{2},\ux\right]=\ux,& \\
& \left[\dD^{2},\ux\right]=-2(1+c)\dD,\qquad&&
\left[\mE+\frac{\delta}{2},\dD^{2}\right]=-2\dD^{2},&\\
& \left[\dD^{2},\ux^{2}\right]=4(1+c)^{2}\left(\mE+\frac{\delta}{2}\right),\qquad&&
\left[\mE+\frac{\delta}{2},\ux^{2}\right]=2\ux^{2},&
\end{alignat*}
where $\delta=1+\frac{\mu-1}{1+c}$.
\end{theorem}

Note that the square of $\dD$ is a~complicated operator, given by
\begin{gather*}
\dD^{2}=-\Delta_{\k}-\left(c\mu\right)r^{-1}\partial_{r}
-\left(c^{2}+2c\right)\partial_{r}^{2}+c r^{-2}\sum_{i}x_{i}T_{i}
-c r^{-2}\sum_{i<j}e_{i}e_{j}(x_{i}T_{j}-x_{j}T_{i}).
\end{gather*}
If $\k=0$, the formula for $\dD^{2}$ simplif\/ies a~bit as now $\sum_{i}x_{i}T_{i}=r
\partial_{r}=\mE$.

\begin{remark}
The operator $\dD=\cD_{\k}+c r^{-2}\ux\mE$ is also considered from a~very dif\/ferent
perspective in~\cite{MR2269930} (in the case $\k=0$), where the eigenfunctions of this operator
are studied.
\end{remark}

Let us now discuss the symmetry of the generators of $\mathfrak{osp}(1|2)$.
First we def\/ine the action of the Pin group on $C^{\infty}(\mR^{m})\otimes\cC l_{m}$ for $s\in
\mathrm{Pin}(m)$ as
\begin{gather*}
\rho(s): \ C^{\infty}(\mR^{m}) \otimes \cC l_{m} \rightarrow C^{\infty}(\mR^{m}) \otimes \cC l_{m},\qquad
 f\otimes b \rightarrow f (p\big(s^{-1}\big)x) \otimes s b.
\end{gather*}
We then have
\begin{proposition}\label{equivarianceG}
Let $s\in\widetilde{\cG}$ and define $\mathrm{sgn}(s):=\mathrm{sgn}(p(s))$.
Then one has
\begin{gather*}
\rho(s) \ux = \mathrm{sgn}(s) \ux \rho(s),\qquad
\rho(s) \dD = \mathrm{sgn}(s) \dD \rho(s).
\end{gather*}
\end{proposition}

\begin{proof}
This follows immediately from the def\/inition of $\rho$ and the $\cG$-equivariance of the Dunkl
operators.
\end{proof}

So up to sign, the Dirac operator $\dD$ is $\widetilde{\cG}$-equivariant.
At this point it is interesting to remark that an algebraic analog of the Dunkl--Dirac operator
$\dD$ for
graded af\/f\/ine Hecke algebras is introduced in~\cite{BCT} with the motivation to prove
a version of Vogan's Conjecture for Dirac cohomology.
The formulation is based on a~uniform geometric parametrization of spin
representations of Weyl groups.
This Dirac operator is an algebraic variant of our family deformation of the dif\/ferential Dirac
operator
for special values of the deformation parameters.
Moreover, it satisf\/ies the same symmetry as in Proposition~\ref{equivarianceG}, see~\cite[Lemma~3.4]{BCT}.

There is a~measure naturally associated with~$\dD$.
Indeed, one has

\begin{proposition}\label{partIntProp}
If $c>-1$, then for suitable differentiable functions $f$ and $g$ one has
\begin{gather*}
\int_{\mR^{m}}\overline{(\dD f)}g h(r)w_{\kappa}(x)dx=\int_{\mR^{m}}\overline{f}
(\dD g)h(r)w_{\k}(x)dx
\end{gather*}
with $h(r)=r^{1-\frac{1+\mu c}{1+c}}$, provided the integrals exist.
\end{proposition}
In this proposition, $\bar{\cdot}$ is the main anti-involution on the Clif\/ford algebra $\cC l_{m}$.

\section[Representation space for the deformation family of the Dunkl-Dirac
operator]{Representation space for the deformation family\\ of the Dunkl--Dirac
operator}\label{section4}

The function space we will work with is $\cL^{2}_{\k,c}(\mR^{m})=L^{2}(\mR^{m},h(r)w_{\k}(x)dx)
\otimes\cC l_{m}$.
This space has the following decomposition
\begin{gather*}
\cL^{2}_{\k,c}(\mR^{m})=L^{2}\big(\mR^{+},r^{\frac{\mu-1}{1+c}}dr\big)\otimes
L^{2}(\mS^{m-1},w_{\k}(\xi)d\sigma(\xi))\otimes\cC l_{m},
\end{gather*}
where on the right-hand side the topological completion of the tensor product is understood and
with $d\sigma(\xi)$ the Lebesgue measure on the sphere $\mS^{m-1}$.
The space $L^{2}(\mS^{m-1},w_{\k}(\xi)d\sigma(\xi))\otimes\cC l_{m}$ can be further decomposed
into Dunkl harmonics and subsequently into Dunkl monogenics.
This leads to
\begin{gather*}
L^{2}\big(\mS^{m-1},w_{\k}(\xi)d\sigma(\xi)\big)\otimes\cC l_{m}=\bigoplus_{\ell=0}^{\infty}
\left(\cM_{\ell}\oplus\ux\cM_{\ell}\right)\big|_{\mS^{m-1}},
\end{gather*}
where $\cM_{\ell}=\ker{\cD_{\k}}\cap\left(\cP_{\ell}\otimes\cC l_{m}\right)$ is the space
of Dunkl monogenics of degree $\ell$, with $\cP_{\ell}$ the space of homogeneous polynomials of
degree $\ell$ (see also~\cite{DBGeg} for more details on Dunkl monogenics).

Using this decomposition, we have obtained in~\cite{H12} a~basis for $\cL^{2}_{\k,c}(\mR^{m})$.
This basis is given by the set $\{\phi_{t,\ell,m}\}$ ($t,\ell\in\mN$ and $m=1,\ldots,\dim
\cM_{\ell}$), def\/ined as
\begin{gather*}
\phi_{2t,\ell,m} = 2^{2t} (1+c)^{2t} t! L_{t}^{\frac{\gamma_{\ell}}{2} -1}(r^{2})r^{\beta_{\ell}}
M_{\ell}^{(m)} e^{-r^{2}/2},\\
\phi_{2t+1,\ell,m} = - 2^{2t+1} (1+c)^{2t+1} t! L_{t}^{\frac{\gamma_{\ell}}{2}}(r^{2}) \ux
r^{\beta_{\ell}} M_{\ell}^{(m)} e^{-r^{2}/2}
\end{gather*}
with $L_{\alpha}^{\beta}$ the Laguerre polynomials and
\begin{gather*}
\beta_{\ell} = - \frac{c}{1+c}\ell,\qquad
\gamma_{\ell} = \frac{2}{1+c}\left( \ell + \frac{\mu-2}{2}\right) + \frac{c+2}{1+c},
\end{gather*}
and where $M_{\ell}^{(m)}$ ($m=1,\ldots,\dim\cM_{\ell}$) forms an orthonormal basis of
$\cM_{\ell}$, i.e.
\begin{gather*}
\left[\int_{\mS^{m-1}}\overline{M_{\ell}^{(m_{1})}}(\xi)M_{\ell}^{(m_{2})}(\xi)w_{\k}(\xi)
d\sigma(\xi)\right]_{0}=\delta_{m_{1}m_{2}}
\end{gather*}
with $[\cdot]_0$ the projection on the scalar part of the Clif\/ford algebra.
The dimension of $\cM_{\ell}$ is given by
\begin{gather*}
\dim_{\mR}{\cM_{\ell}} = \dim_{\mR}{\cC l_{m}} \dim_{\mR}{\cP_{\ell}\big(\mR^{m-1}\big)}
=2^{m} \frac{(\ell +m-2)!}{\ell! (m-2)!}
\end{gather*}
with $\cP_{\ell}\big(\mR^{m-1}\big)$ the space of homogeneous polynomials of degree $\ell$ in $m-1$
variables (see~\cite{MR1169463}).

Using formula (4.10) in~\cite{H12} and the proof of Theorem 3 in \cite{H12}, one obtains the
following formulae for the action of $\dD$ and $\ux$ on the generalized Laguerre functions
\begin{gather}
2 \dD \phi_{t,\ell,m}\! = \phi_{t+1,\ell,m}\! + C(t,\ell) \phi_{t-1,\ell,m},\qquad
-2(1+c) \ux \phi_{t,\ell,m}\! = \phi_{t+1,\ell,m}\! - C(t,\ell) \phi_{t-1,\ell,m}\!\!\!\!\label{ActionOnBasis}
\end{gather}
with
\begin{gather*}
C(2t, \ell)= 4 (1+c)^{2} t,\qquad
C(2t+1, \ell)= 2 (1+c)^{2}(\gamma_{\ell}+ 2t).
\end{gather*}
These formulae determine the action of $\mathfrak{osp}(1|2)$ on $\cL^{2}_{\k,c}(\mR^{m})$.
Recall also that the action of $\widetilde{\cG}$ on $\cL^{2}_{\k,c}(\mR^{m})$ is given by $\rho$
(see Section~\ref{section3}).

Subsequently, we can def\/ine a~creation and annihilation operator in this setting by
\begin{gather}\label{CreaAnn}
A^{+} = \dD - (1+c) \ux,\qquad
A^{-} = \dD + (1+c) \ux
\end{gather}
satisfying
\begin{gather*}
A^{+} \phi_{t,\ell,m} = \phi_{t+1,\ell,m},\qquad
A^{-} \phi_{t,\ell,m} = C(t, \ell) \phi_{t-1,\ell,m}.
\end{gather*}

Now we introduce the following inner product
\begin{gather*}
\langle f,g\rangle=\left[\int_{\mR^{m}}\overline{f^{c}}g h(r)w_{\k}(x)dx\right]_{0},
\end{gather*}
where $h(r)$ is the measure associated to $\dD$ (see Proposition~\ref{partIntProp}) and $f^{c}$ is
the complex conjugate of $f$.
It is easy to check that this inner product satisf\/ies
\begin{gather}\label{adjointsIP}
\langle \dD f, g \rangle = \la f, \dD g\ra,\qquad
\la \ux f, g \ra = -\la f, \ux g \ra.
\end{gather}
The related norm is def\/ined by $||f||^{2}=\la f,f\ra$.

\begin{theorem}
We have
\begin{gather*}
\langle\phi_{t_{1},\ell_{1},m_{1}},\phi_{t_{2},\ell_{2},m_{2}}\rangle=c(t_{1},\ell_{1})
\delta_{t_{1}t_{2}}\delta_{\ell_{1}\ell_{2}}\delta_{m_{1}m_{2}},
\end{gather*}
where $c(t,\ell)$ is a~constant depending on $t$ and $\ell$.
\end{theorem}

The functions $\phi_{t,\ell,m}$ are eigenfunctions of the Hamiltonian of a~generalized harmonic
oscillator.
\begin{theorem}
\label{HarmOsc}
The functions $\phi_{t,\ell,m}$ satisfy the following second-order PDE
\begin{gather*}
\left(\dD^{2}-(1+c)^{2}\ux^{2}\right)\phi_{t,\ell,m}=(1+c)^{2}(\gamma_{\ell}+2t)
\phi_{t,\ell,m}.
\end{gather*}
\end{theorem}

\begin{proof}
This follows immediately from the formula~\eqref{ActionOnBasis}.
\end{proof}

Theorem~\ref{HarmOsc} combined with the def\/inition of $A^{+}$, $A^{-}$ in~\eqref{CreaAnn} allows us to
decompose the space $\cL^{2}_{\k,c}(\mR^{m})$ under the action of $\mathfrak{osp}(1|2)$.
Clearly the odd elements $A^{+}$ and $A^{-}$ generate
$\mathfrak{osp}(1|2)$ as they are linear combinations of $\dD$ and $\ux$.
Moreover, they act between two basis vectors
$\{\phi_{t,\ell,m}\}$ of $\cL^{2}_{\k,c}(\mR^{m})$, so it is suf\/f\/icient to consider vectors in an
irreducible representation of $\mathfrak{osp}(1|2)$ inside the functional space.
This is achieved as follows~-- for f\/ixed $\ell$ and $m$ each vector $\phi_{0,\ell,m}$ generates the
irreducible representation
\begin{gather*}
\xymatrix{
\phi_{0,\ell,m}\ar@<1ex>[r]^{A^+}\ar@(dl,dr)_{L}&\phi_{1,\ell,m}\ar@(dl,dr)_{L}
\ar@<1ex>[r]^{A^+}\ar@<1ex>[l]^{A^-}&\phi_{2,\ell,m}\ar@(dl,dr)_{L}\ar@<1ex>[r]^{A^+}
\ar@<1ex>[l]^{A^-}&\phi_{3,\ell,m}\ar@(dl,dr)_{L}\ar@<1ex>[r]^{A^+}\ar@<1ex>[l]^{A^-}
&\phi_{4,\ell,m}\ar@(dl,dr)_{L}\ar@<1ex>[r]^{A^+}\ar@<1ex>[l]^{A^-}&\ar@<1ex>[l]^{A^-}\ldots
}
\end{gather*}
where
\begin{gather*}
L=\frac{1}{2}\{A^{+},A^{-}\}=\dD^{2}-(1+c)^{2}\ux^{2}
\end{gather*}
with the action given in Theorem~\ref{HarmOsc}.
In fact this highest weight representation is labeled by~$\ell$ only and we will denote it
$\pi(\ell)$.
In conclusion, we obtain the decomposition of our functional space $\cL^{2}_{\k,c}(\mR^{m})$
into a~discrete direct sum of highest weight (inf\/inite-dimensional)
Harish-Chandra modules for~$\mathfrak{osp}(1|2)$:
\begin{gather*}
\cL^{2}_{\k,c}(\mR^{m})=\bigoplus_{\ell=0}^{\infty}\pi(\ell)\otimes\cM_{\ell}.
\end{gather*}

These results should be compared with Theorem 3.19 and Section~3.6 in~\cite{Orsted2} (where one
uses~$\mathfrak{sl}_{2}$ instead of $\mathfrak{osp}(1|2)$).
Also notice that the claim should be understood as an assertion on the deformation of the Howe dual
pair for $\mathfrak{osp}(1|2)$ inside the Clif\/ford--Weyl algebra on $\mR^{m}$ acting on a~f\/ixed vector
space $\cL^{2}_{\k,c}(\mR^{m})$.

In particular, we have the following result.
Recall that an operator $T$ is essentially selfadjoint on a~Hilbert space $H$ if $T$ is a
symmetric operator with a~dense domain $D(T)\subset H$ such that for a~complete orthogonal set
$\{f_n\}_n$ in $H$ with $f_n\in D(H)$, there exist $\{\mu_n\}_n$ solving $Tf_n=\mu_nf$ for all
$n\in\mN$.
\begin{proposition}
Let $c>-1$ and $\kappa>0$.
The operator $L$ acting on
$\cL^{2}_{\k,c}(\mR^{m})$ is essentially selfadjoint $($i.e.\
symmetric and its closure is a~selfadjoint operator$)$.
Moreover, $L$ has no continuous spectrum and its
discrete spectrum is given by
\begin{gather*}
\operatorname{Spec}(L)=\{2(1+c)\ell+2(1+c)^{2}t+(1+c)(\mu+c)\,|\,\ell,t\in\mN\}.
\end{gather*}
\end{proposition}

Using Theorem~\ref{HarmOsc} we can now def\/ine the holomorphic semigroup for the deformed Dirac
operator by
\begin{gather*}
\cF_{\dD}^{\omega}=e^{\o\left(\frac{1}{2}+\frac{\mu-1}{2(1+c)}\right)}e^{\frac{-\o}{
2(1+c)^{2}}\left(\dD^{2}-(1+c)^{2}\ux^{2}\right)}.
\end{gather*}
Here, $\omega$ takes values in the right half-plane of $\mC$.
The special boundary value $\o=i\pi/2$ corresponds to the Fourier transform.
In that case, we will use the notation $\cF_{\dD}$.
The functions $\phi_{t,\ell,m}$ are eigenfunctions of $\cF_{\dD}^{\omega}$ satisfying
\begin{gather}
\label{eigvalsFourier}
\cF_{\dD}^{\o}(\phi_{t,\ell,m})=e^{-\o t}e^{-\frac{\o\ell}{(1+c)}}\phi_{t,\ell,m}.
\end{gather}
Note that in the special case $\kappa=0$, $c=0$ the operator $\cF_{\dD}^{\omega}$ reduces to the
classical Hermite semigroup (see, e.g.,~\cite{MR0974332}).

\begin{remark}
One can also consider more general deformations of the Dirac operator, by adding suitable odd
powers of $\Gamma=-\ux\cD_{\k}-\mE$ to $\dD$ as follows
\begin{gather*}
\dD=\cD_{\k}+c r^{-2}\ux\mE+\sum_{j=0}^{\ell}c_{j}r^{-1}\left(\Gamma-\frac{\mu-1}{2}
\right)^{2j+1},\qquad c_{j}\in\mR.
\end{gather*}
This does not alter the $\mathfrak{osp}(1|2)$ relations, as $\Gamma-\frac{\mu-1}{2}$
anti-commutes with $\ux$ and has the correct homogeneity.
In particular, $\Gamma-\frac{\mu-1}{2}$ can be seen as the square root of the Casimir of
$\mathfrak{osp}(1|2)$, see~\cite[Example~2 in Section~2.5]{MR1773773}.
\end{remark}

In the sequel of the paper, we will always assume $\kappa=0$ or in other words, we do not consider
the Dunkl deformation.
This is to simplify the notation of the results.
Most statements can be generalized to the Dunkl case by a~suitable composition with the Dunkl
intertwining operator~$V_{\k}$, except the results obtained in Section~\ref{section8}.

Recall that for $\kappa=0$, the Dunkl--Dirac operator $\cD_{\k}$ reduces to the orthogonal Dirac
operator $\upx=\sum\limits_{i=1}^{m}e_{i}\partial_{x_{i}}$ and the Dunkl dimension~$\mu$ to the ordinary
dimension~$m$.

\section{Reproducing kernels}\label{section5}

In this section we determine the reproducing kernels for $\cM_{k}$ and $\ux\cM_{k}$.
We start with an auxiliary Lemma, which can be thought of as a~Clif\/ford analogue
of the Funk--Hecke transform.
We def\/ine the wedge product of two vectors as
\begin{gather*}
\ux\wedge\uy:=\sum_{j<k}e_{j}e_{k}(x_{j}y_{k}-x_{k}y_{j}).
\end{gather*}

\begin{lemma}
\label{auxintegrals}
Put $\ux=r\ux'$ and $\uy=s\uy'$ with $\ux',\uy'\in\mS^{m-1}$.
Furthermore, put $\lambda=(m-2)/2$ and $\sigma_m=2\pi^{m/2}/\Gamma(m/2)$.
Then one has, with $M_{l}\in\cM_{\ell}$
\begin{gather*}
\int_{\mS^{m-1}} C_{k}^{\lambda}(\la\ux', \uy' \ra) M_{\ell}(\ux')d \s (\ux') = \sigma_m
\frac{\l}{\l+k} \delta_{k, \ell} M_{\ell}(\uy'),\\
\int_{\mS^{m-1}} C_{k}^{\lambda}(\la\ux', \uy' \ra) \ux' M_{\ell}(\ux')d \s (\ux') = \sigma_m
\frac{\l}{\l+k} \delta_{k, \ell+1} \uy' M_{\ell}(\uy'),\\
\int_{\mS^{m-1}} (\ux' \wedge \uy' ) C_{k-1}^{\lambda+1}(\la\ux', \uy' \ra) M_{\ell}(\ux')d \s
(\ux') = - \sigma_m \frac{k}{2(\l+k)} \delta_{k, \ell} M_{\ell}(\uy'),\\
\int_{\mS^{m-1}} (\ux' \wedge \uy' ) C_{k-1}^{\lambda+1}(\la\ux', \uy' \ra) \ux' M_{\ell}(\ux')d
\s (\ux') = \sigma_m \frac{k+2\l}{2(\l+k)} \delta_{k, \ell+1} \uy' M_{\ell}(\uy'),
\end{gather*}
where $C_{k}^{\lambda}(\la\ux',\uy'\ra)$ is the $k$-th Gegenbauer polynomial in the variable
$\la\ux',\uy'\ra$.
\end{lemma}

\begin{proof}
The f\/irst integral is trivial: $M_{\ell}$ is a~spherical harmonic of degree $\ell$ and
$C_{k}^{\lambda}(\la\ux',\uy'\ra)$ is the reproducing kernel for spherical harmonics of degree
$k$ (see, e.g.,~\cite{MR1827871}).
The second integral immediately follows, because $\ux'M_{\ell}(\ux')\in\cH_{\ell+1}$.

The other two integrals are a~bit more complicated.
We show how to obtain the last one.
First rewrite
$(\ux'\wedge\uy')\ux'=\uy'-\la\ux',\uy'\ra\ux'$.
The f\/irst term then follows from the f\/irst integral.
For the second term, we use the recursive property of Gegenbauer polynomials:
\begin{gather*}
w C^{\l+1}_{n-1}(w)=\frac{n}{2(n+\l)}C^{\l+1}_{n}(w)+\frac{n+2\l}{2(n+\l)}
C^{\l+1}_{n-2}(w).
\end{gather*}
The result then follows by collecting everything.
\end{proof}

We can use this lemma to determine the reproducing kernels.
This is the subject of the following proposition.

\begin{proposition}
\label{reprkernels}
For $k\in\mN^{*}$ put
\begin{gather*}
P_{k}(\ux', \uy') = \frac{k+2\l}{2\l} C_{k}^{\l}(\la\ux', \uy' \ra) -(\ux' \wedge \uy' )
C_{k-1}^{\lambda+1}(\la\ux', \uy' \ra),\\
Q_{k-1}(\ux', \uy') = \frac{k}{2\l} C_{k}^{\l}(\la\ux', \uy' \ra) +(\ux' \wedge \uy' )
C_{k-1}^{\lambda+1}(\la\ux', \uy' \ra)
\end{gather*}
with $P_{0}(\ux',\uy')=C_{0}^{\l}(0)=1$.
Then
\begin{gather*}
\int_{\mS^{m-1}} P_{k}(\ux', \uy') M_{\ell}(\ux')d \s (\ux') = \sigma_m \delta_{k, \ell}
M_{\ell}(\uy'),\\
\int_{\mS^{m-1}} P_{k}(\ux', \uy') \ux'M_{\ell}(\ux')d \s (\ux') = 0
\end{gather*}
and
\begin{gather*}
\int_{\mS^{m-1}} Q_{k-1}(\ux', \uy') M_{\ell}(\ux')d \s (\ux') = 0,\\
\int_{\mS^{m-1}} Q_{k-1}(\ux', \uy') \ux'M_{\ell}(\ux')d \s (\ux') = \sigma_m \delta_{k, \ell+1}
\uy'M_{\ell}(\uy').
\end{gather*}
\end{proposition}

\begin{proof}
This follows immediately from Lemma~\ref{auxintegrals}.
\end{proof}

\begin{remark}
Note that, as expected, $P_{k}(\ux',\uy')+Q_{k-1}(\ux',\uy')=
\frac{\l+k}{\l}C_{k}^{\l}(\la\ux',\uy'\ra)$, which is the reproducing kernel for the space of
spherical harmonics of degree $k$.
\end{remark}

\begin{remark}
When the dimension $m=2$ and hence $\lambda=0$, the reproducing kernel is still well-def\/ined by
using the well-known relation~\cite[(4.7.8)]{Sz}
\begin{gather*}
\lim_{\lambda\rightarrow0}\lambda^{-1}C_k^\lambda(w)=(2/k)\cos k\theta,\qquad w=\cos
\theta,\qquad k\geq1.
\end{gather*}
\end{remark}

We will also need the following lemma.

\begin{lemma}
\label{ReprKernOrth}
The reproducing kernels satisfy the following properties, for all $k,l\in\mN$:
\begin{gather*}
\int_{\mS^{m-1}} P_{k}(\uy', \ux') P_{\ell}(\uz', \uy') d \s (\uy') = \sigma_m \delta_{k \ell
}P_{\ell}( \uz', \ux'),\\
\int_{\mS^{m-1}} P_{k}( \uy', \ux') Q_{\ell}( \uz', \uy') d \s (\uy') = 0,\\
\int_{\mS^{m-1}} Q_{k}( \uy', \ux') Q_{\ell}(\uz', \uy') d \s (\uy') = \sigma_m \delta_{k \ell }
Q_{\ell}(\uz', \ux').
\end{gather*}
\end{lemma}

\begin{proof}
This follows immediately using Lemma 7.6 and 7.10 from~\cite{DBXu}.
\end{proof}

\begin{remark}
Mind the order of the variables in the previous lemma.
The kernels $P_{k}(\ux',\uy')$ and $Q_{k}(\ux',\uy')$ are not symmetric.
\end{remark}

\section{The series representation of the holomorphic semigroup}\label{section6}

The aim of the present section is to investigate basic properties of the holomorphic semigroup
def\/ined by
\begin{gather*}
\cF_{\dD}^{\o}=e^{\o\left(\frac{1}{2}+\frac{\mu-1}{2(1+c)}\right)}e^{\frac{-\o}{
2(1+c)^{2}}\left(\dD^{2}-(1+c)^{2}\ux^{2}\right)},\qquad\operatorname{Re}\o\geq0,
\end{gather*}
acting on the space $\cL^{2}_{0,c}(\mR^{m})$.
We start with the following general statement.

\begin{theorem}\label{PropSemigroup}
Suppose $c>-1$.
Then
\begin{enumerate}\itemsep=0pt
\item[$1.$]
For any $t,\ell\in\mN$ and $m\in\{1,\ldots,\dim\cM_{\ell}\}$, the function
$\phi_{t,\ell,m}$ is an eigenfunction of the
operator $\cF_{\dD}^{\o}$:
\begin{gather*}
\cF_{\dD}^{\o}(\phi_{t,\ell,m})=e^{-\o t}e^{-\frac{\o\ell}{(1+c)}}\phi_{t,\ell,m}.
\end{gather*}
\item[$2.$] $\cF_{\dD}^{\o}$ is a~continuous operator on $\cL^{2}_{0,c}(\mR^{m})$ for all $\o$ with $\operatorname{Re}
\o\geq0$, in particular
\begin{gather*}
||\cF_{\dD}^{\o}(f)||\leq||f||
\end{gather*}
for all $f\in\cL^{2}_{0,c}(\mR^{m})$.
\item[$3.$]
If $\operatorname{Re}\o>0$, then $\cF_{\dD}^{\o}$ is a~Hilbert--Schmidt operator on $\cL^{2}_{0,c}(\mR^{m})$.
\item[$4.$]
If $\operatorname{Re}\o=0$, then $\cF_{\dD}^{\o}$ is a~unitary operator on $\cL^{2}_{0,c}(\mR^{m})$.
\end{enumerate}
\end{theorem}

\begin{proof}
$(1)$ is an immediate consequence of Theorem~\ref{HarmOsc}.
For $(2)$, let $f$ be an element in $\cL^{2}_{0,c}(\mR^{m})$ and expand it with respect to the
(normalized) basis $\{\phi_{t,\ell,m}\}$ as
\begin{gather*}
f=\sum_{t,\ell,m}a_{t,\ell,m}\phi_{t,\ell,m}.
\end{gather*}
Then one has, using orthogonality,
\begin{gather*}
||\cF_{\dD}^{\o}(f)||^2 = \sum_{t, \ell, m} |a_{t,\ell,m}|^{2} e^{-2 (\operatorname{Re}{\o}) t} e^{-\frac{2(\operatorname{Re}{\o}) \ell}{(1+c)}}
\leq \sum_{t, \ell, m} |a_{t,\ell,m}|^{2} = ||f||^2
\end{gather*}
because $\operatorname{Re}{\o}\geq0$.

As for $(3)$, we have to show that the Hilbert--Schmidt norm is f\/inite.
We compute
\begin{gather*}
||\cF_{\dD}^{\o}||_{\rm HS}^{2}
= \sum_{t, k, \ell} ||\cF_{\dD}^{\o}(\phi_{t,k,\ell})||^{2}
= \sum_{t, k,\ell} e^{-2 (\operatorname{Re}{\o}) t} e^{-\frac{2 (\operatorname{Re}{\o}) k}{(1+c)}}
= \sum_{t, k} e^{-2 (\operatorname{Re}{\o}) t} e^{-\frac{2 (\operatorname{Re}{\o}) k}{(1+c)}} \dim_{\mR}{\cM_{k}}\\
{} = \sum_{t, k} e^{-2 (\operatorname{Re}{\o}) t} e^{-\frac{2 (\operatorname{Re}{\o}) k}{(1+c)}} \frac{(k + m - 2)!}{k!(m - 2)!}2^{m}
= \sum_{t} e^{-2 (\operatorname{Re}{\o}) t} \sum_{k} e^{-\frac{2 (\operatorname{Re}{\o}) k}{(1+c)}} \frac{(k + m - 2)!}{k!(m - 2)!}2^{m}.
\end{gather*}
Using the ratio test, we see that these series are convergent for $\operatorname{Re}{\o}>0$.

$(4)$ follows immediately, because when $\operatorname{Re}{\o}=0$ the eigenvalues all have unit norm.
\end{proof}

We have already observed that $\cF_{\dD}^{\o}$ is a~Hilbert--Schmidt operator
for $\operatorname{Re}\o>0$ and a~unitary operator for $\operatorname{Re}\o=0$.
The Schwartz kernel theorem implies
that $\cF_{\dD}^{\o}$ can be expressed by a~distribution kernel
$K(x,y;\omega)$, so
\begin{gather*}
\big(\cF_{\dD}^{\o}f\big)(y)=\int_{\mR^m}K(x,y;\omega)f(x)h(r_{x})dx,
\end{gather*}
and $K(x,y;\omega)h(r_{x})$ is a~tempered distribution
on $\mR^m\times\mR^m$.

\subsection[The case $\operatorname{Re}\o>0$]{The case $\boldsymbol{\operatorname{Re}\o>0}$}

Using the reproducing kernels of Section~\ref{section5}, we can now make a~reasonable ansatz for
the kernel of the full holomorphic semigroup.
We want to write this semigroup as
\begin{gather*}
\cF_{0,c}^{\o}(f)(y)=\sigma_{m}^{-1}\int_{\mR^{m}}K(x,y;\omega)f(x)h(r_{x})dx
\end{gather*}
with $K(x,y;\omega)=K_{0}(x,y;\omega)+K_{1}(x,y;\omega)$ and
\begin{gather}
K_{0} = e^{-\frac{\coth{\o}}{2} (r^{2}+ s^{2})}\sum_{k=0}^{+\infty} \alpha_{k} z^{\frac{k}{1+c}}
\widetilde{J}_{\frac{\gamma_{k}}{2} -1}\left(\frac{iz}{\sinh{\o}}\right) P_{k}(\ux', \uy'),\nonumber\\
K_{1} = e^{-\frac{\coth{\o}}{2} (r^{2}+ s^{2})} \sum_{k=0}^{+\infty} \beta_{k} z^{1+\frac{k}{1+c}}
\widetilde{J}_{\frac{\gamma_{k}}{2}}\left(\frac{iz}{\sinh{\o}}\right) Q_{k}(\ux', \uy').\label{seriessemigroup}
\end{gather}
Here $r=|\ux|$, $s=|\uy|$ and $z=|\ux||\uy|$.
We also used the notation $\widetilde{J}_{\nu}(t)=(t/2)^{-\nu}J_{\nu}(t)$.
Now we determine the complex constants $\{\alpha_{k}\}$ and $\{\beta_{k}\}$ such that this
integral transform coincides with
\begin{gather*}
\cF_{\dD}^{\o}=e^{\o\left(\frac{1}{2}+\frac{\mu-1}{2(1+c)}\right)}e^{\frac{-\o}{
2(1+c)^{2}}\left(\dD^{2}-(1+c)^{2}\ux^{2}\right)}
\end{gather*}
on the basis $\{\phi_{t,\ell,m}\}$.

We calculate
\begin{gather*}
\sigma_{m}^{-1}\int_{\mR^{m}} K_{0}(x,y;\omega) \phi_{2t,\ell,m} (x)dx = \alpha_{\ell}
M_{\ell}^{(m)}(\uy') e^{-\frac{\coth{\o}}{2}s^{2}}s^{\frac{\ell}{1+c}}\left( \frac{i s}{2
\sinh{\o}}\right)^{-\gamma_{\ell}/2+1}\\
{}
\times \int_{0}^{+\infty} r^{\gamma_{\ell}/2} e^{-\frac{(\coth{\o}+1)}{2}r^{2}}
J_{\gamma_{\ell}/2 -1} \left( \frac{i r s}{\sinh{\o}}\right)L_{t}^{\gamma_{\ell}/2-1}(r^{2}) dr
=\alpha_{\ell} \frac{e^{-2\o t} 2^{\gamma_{\ell}/2-1}}{(\coth{\o}+1)^{\gamma_{\ell}/2}}
\phi_{2t,\ell,m}(y),
\end{gather*}
where we used the identity (see~\cite[Corollary 4.6]{Orsted2})
\begin{gather*}
2\int_{0}^{+\infty}r^{\alpha+1}J_{\a}(r\beta)L_{j}^{\a}(r^{2})e^{-\d r^{2}}dr=
\frac{(\d-1)^{j}\beta^{\alpha}}{2^{\alpha}\d^{\alpha+j+1}}L_{j}^{\a}\left(\frac{\beta^{2}}{4
\d(1-\d)}\right)e^{-\frac{\beta^{2}}{4\d}}.
\end{gather*}
Similarly, we f\/ind
\begin{gather*}
\sigma_{m}^{-1} \int_{\mR^{m}} K_{0}(x,y;\omega) \phi_{2t+1,\ell,m}(x) dx = 0,\qquad
\sigma_{m}^{-1} \int_{\mR^{m}} K_{1}(x,y;\omega) \phi_{2t,\ell,m} (x)dx = 0,\\
\sigma_{m}^{-1} \int_{\mR^{m}} K_{1}(x,y;\omega) \phi_{2t+1,\ell,m} (x) dx = \beta_{\ell}
\frac{e^{-2\o t} 2^{\gamma_{\ell}/2}}{(\coth{\o}+1)^{\gamma_{\ell}/2+1}} \phi_{2t+1,\ell,m}(y).
\end{gather*}
Hence we obtain by comparison with~\eqref{eigvalsFourier}
\begin{gather*}
\alpha_{\ell} = e^{-\frac{\o \ell}{(1+c)}}
\frac{(\coth{\o}+1)^{\gamma_{\ell}/2}}{2^{\gamma_{\ell}/2-1}} = 2 e^{\frac{\o \d}{2}} (2
\sinh{\o})^{-\gamma_{\ell}/2},\qquad
\beta_{\ell} = \frac{ \alpha_{\ell}}{2\sinh{\o}}.
\end{gather*}

We summarize our results in the following theorem.

\begin{theorem}
\label{seriesholom}
Let $\operatorname{Re}\o>0$ and $c>-1$.
Put
\begin{gather*}K(x,y;\omega)=e^{-\frac{\coth{\o}}{2}(r^{2}+s^{2})}\left(A(z,w)+\ux\wedge\uy
B(z,w)\right)
\end{gather*} with
\begin{gather*}
A(z,w) = \sum_{k=0}^{+\infty} \left(\alpha_{k} \frac{k+2\l}{2\l} z^{\frac{k}{1+c}}
\widetilde{J}_{\frac{\gamma_{k}}{2}-1}\left(\frac{iz}{\sinh{\o}}\right) +\frac{
\alpha_{k-1}}{4\sinh{\o}} \frac{k}{\l} z^{\frac{k+c}{1+c}}
\widetilde{J}_{\frac{\gamma_{k-1}}{2}}\left(\frac{iz}{\sinh{\o}}\right)\right) C_{k}^{\l}(w),\\
B(z,w) = \sum_{k=1}^{+\infty} \left(-\alpha_{k} z^{\frac{k}{1+c}-1}
\widetilde{J}_{\frac{\gamma_{k}}{2}-1}\left(\frac{iz}{\sinh{\o}}\right) +
\frac{\alpha_{k-1}}{2\sinh{\o}} z^{\frac{k+c}{1+c}-1}
\widetilde{J}_{\frac{\gamma_{k-1}}{2}}\left(\frac{iz}{\sinh{\o}}\right)\right) C_{k-1}^{\l+1}(w)
\end{gather*}
for $z=|\ux||\uy|$, $w=\la\ux,\uy\ra/z$, $\alpha_{-1}=0$ and
$
\alpha_{k}=2e^{\frac{\o\d}{2}}(2\sinh{\o})^{-\gamma_{k}/2}.
$

Then these series are convergent and the integral transform defined on $\cL^{2}_{0,c}(\mR^{m})$ by
\begin{gather*}
\cF_{0,c}^{\o}(f)(y)=\sigma_{m}^{-1}\int_{\mR^{m}}K(x,y;\omega)f(x)h(r_{x})dx
\end{gather*}
coincides with the operator $\cF_{\dD}^{\o}=e^{\o\left(\frac{1}{2}+\frac{\mu-1}{2(1+c)}
\right)}e^{\frac{-\o}{2(1+c)^{2}}\left(\dD^{2}-(1+c)^{2}\ux^{2}\right)}$ on the basis $\{
\phi_{t,\ell,m}\}$.
\end{theorem}

\begin{proof}
We have already shown that the integral transform coincides with the operator $\cF_{\dD}^{\o}=
e^{\o\left(\frac{1}{2}+\frac{\mu-1}{2(1+c)}\right)}e^{\frac{-\o}{2(1+c)^{2}}\left(\dD^{2}-
(1+c)^{2}\ux^{2}\right)}$ on the basis $\{\phi_{t,\ell,m}\}$.
So we only have to show that the series are convergent.
We do this for the term
\begin{gather*}
\sum_{k=0}^{+\infty} (2 \sinh{\o})^{-\gamma_{k}/2} \frac{k+2\l}{2\l} z^{\frac{k}{1+c}}
\widetilde{J}_{\frac{\gamma_{k}}{2}-1}\left(\frac{iz}{\sinh{\o}}\right) C_{k}^{\l}(w)\\
\qquad
{}= (2 \sinh{\o})^{-\d/2} \sum_{k=0}^{+\infty} \frac{k+2\l}{2\l} \left(\frac{z}{2
\sinh{\o}}\right)^{\frac{k}{1+c}}
\widetilde{J}_{\frac{\gamma_{k}}{2}-1}\left(\frac{iz}{\sinh{\o}}\right) C_{k}^{\l}(w),
\end{gather*}
the other ones are treated in a~similar fashion.
We obtain
\begin{gather*}
\left| \sum_{k=0}^{+\infty} \frac{k+2\l}{2\l} \left(\frac{z}{2 \sinh{\o}}\right)^{\frac{k}{1+c}}
\widetilde{J}_{\frac{\gamma_{k}}{2}-1}\left(\frac{iz}{\sinh{\o}}\right) C_{k}^{\l}(w) \right| \\
\qquad
{}\leq \frac{B(\l)}{2}e^{\left|\operatorname{Im}{\frac{iz}{\sinh{\o}}}\right|}\sum_{k=0}^{+\infty}(k+2\l)\left|\frac{z}{2
\sinh{\o}}\right|^{\frac{k}{1+c}} \frac{1}{\Gamma(\gamma_{k}/2)} k^{2\l-1}
\end{gather*}
using formula~\eqref{estimateGeg} and \eqref{estimateBessel}.
As the term $\Gamma(\gamma_{k}/2)$ is dominant, the series clearly converges.
\end{proof}

\subsection[The case $\operatorname{Re}\o=0$]{The case $\boldsymbol{\operatorname{Re}\o=0}$}

In this case, we have the following theorem.

\begin{theorem}
Let $c>-1$.
Then for $\o=i\eta$ with $\eta\not\in\pi\mZ$, we put
\begin{gather*}K(x,y;i\eta)=e^{i\frac{\cot{\eta}}{2}(r^{2}+s^{2})}\left(A(z,w)+\ux\wedge\uy
B(z,w)\right)
\end{gather*}
with{\samepage
\begin{gather*}
A(z,w) = \sum_{k=0}^{+\infty} \left(\alpha_{k} \frac{k+2\l}{2\l} z^{\frac{k}{1+c}}
\widetilde{J}_{\frac{\gamma_{k}}{2}-1}\left(\frac{z}{\sin{\eta}}\right) +\frac{
\alpha_{k-1}}{4i\sin{\eta}} \frac{k}{\l} z^{\frac{k+c}{1+c}}
\widetilde{J}_{\frac{\gamma_{k-1}}{2}}\left(\frac{z}{\sin{\eta}}\right)\right) C_{k}^{\l}(w),\\
B(z,w) = \sum_{k=1}^{+\infty} \left(-\alpha_{k} z^{\frac{k}{1+c}-1}
\widetilde{J}_{\frac{\gamma_{k}}{2}-1}\left(\frac{z}{\sin{\eta}}\right) +
\frac{\alpha_{k-1}}{2i\sin{\eta}} z^{\frac{k+c}{1+c}-1}
\widetilde{J}_{\frac{\gamma_{k-1}}{2}}\left(\frac{z}{\sin{\eta}}\right)\right) C_{k-1}^{\l+1}(w)
\end{gather*}
for $z=|\ux||\uy|$, $w=\la\ux,\uy\ra/z$, $\alpha_{-1}=0$ and
$
\alpha_{k}=2e^{\frac{i\eta\d}{2}}(2i\sin{\eta})^{-\gamma_{k}/2}$.}

These series are convergent and the unitary integral transform defined in distributional sense on
$\cL^{2}_{0,c}(\mR^{m})$ by
\begin{gather*}
\cF_{0,c}^{i\eta}(f)(y)=\sigma_{m}^{-1}\int_{\mR^{m}}K(x,y;i\eta)f(x)h(r_{x})dx
\end{gather*}
coincides with the operator $\cF_{\dD}^{i\eta}=e^{i\eta\left(\frac{1}{2}+\frac{\mu-1}{2(1+c)}
\right)}e^{\frac{-i\eta}{2(1+c)^{2}}\left(\dD^{2}-(1+c)^{2}\ux^{2}\right)}$ on the basis $\{
\phi_{t,\ell,m}\}$.
\end{theorem}

\begin{proof}
This follows by taking the limit $\o\rightarrow i\eta$ in Theorem~\ref{seriesholom}.
\end{proof}

\section{The series representation of the Fourier transform}\label{section7}

The Fourier transform is the very special case of the holomorphic semigroup, evaluated at $\o=i
\pi/2$.
In this case, the kernel $K(x,y)=K(x,y;i\pi/2)$ is given by the following theorem.

\begin{theorem}
\label{seriesrepr}
Put $K(x,y)=A(z,w)+\ux\wedge\uy B(z,w)$ with
\begin{gather*}
A(z,w) = \sum_{k=0}^{+\infty} z^{-\frac{\delta-2}{2}} \left(\alpha_{k} \frac{k+2\l}{2\l}
J_{\frac{\gamma_{k}}{2}-1}(z) -i \alpha_{k-1} \frac{k}{2\l} J_{\frac{\gamma_{k-1}}{2}}(z)\right)
C_{k}^{\l}(w),\\
B(z,w) = \sum_{k=1}^{+\infty} z^{-\frac{\delta}{2}} \left(-\alpha_{k}
J_{\frac{\gamma_{k}}{2}-1}(z) -i \alpha_{k-1} J_{\frac{\gamma_{k-1}}{2}}(z)\right)
C_{k-1}^{\l+1}(w)
\end{gather*}
and $z=|\ux||\uy|$, $w=\la\ux,\uy\ra/z$, $\alpha_{-1}=0$ and $\alpha_{k}=e^{-\frac{i\pi
k}{2(1+c)}}$.
These series are convergent and the integral transform defined in distributional sense on
$\cL^{2}_{0,c}(\mR^{m})$ by
\begin{gather*}
\cF_{0,c}(f)(y)=\sigma_{m}^{-1}\int_{\mR^{m}}K(x,y)f(x)h(r_{x})dx
\end{gather*}
coincides with the operator $\cF_{\dD}=e^{i\frac{\pi}{2}\left(\frac{1}{2}+\frac{\mu-
1}{2(1+c)}\right)}e^{\frac{-i\pi}{4(1+c)^{2}}\left(\dD^{2}-(1+c)^{2}\ux^{2}\right)}$
on the basis $\{\phi_{t,\ell,m}\}$.
\end{theorem}

\begin{proof}
Using the well-known identity (see~\cite[Exercise 21, p.~371]{Sz})
\begin{gather*}
\int_{0}^{+\infty}r^{\alpha+1}J_{\a}(rs)L_{j}^{\a}\big(r^{2}\big)e^{-r^{2}/2}dr=(-1)^{j}s^{\a}
L_{j}^{\a}(s^{2})e^{-s^{2}/2}
\end{gather*}
we can prove in the same way as leading to Theorem~\ref{seriesholom} that the integral transform
$\cF_{0,c}$ coincides with
\begin{gather*}
\cF_{\dD}=e^{i\frac{\pi}{2}\left(\frac{1}{2}+\frac{\mu-1}{2(1+c)}\right)}e^{\frac{-i\pi}{
4(1+c)^{2}}\left(\dD^{2}-(1+c)^{2}\ux^{2}\right)}
\end{gather*}
on the basis $\phi_{t,\ell,m}$.
The theorem also follows as a~special case of Theorem~\ref{seriesholom}, taking the limit $\o
\rightarrow i\pi/2$.
\end{proof}

\begin{remark}
One can also def\/ine an analogue of the Schwartz space of rapidly decreasing functions in this
context.
Let $L=\dD^{2}-(1+c)^{2}\ux^{2}$ and denote by $D(L)$ the domain of~$L$ in~$\cL^{2}_{0,c}(\mR^{m})$.
Then the Schwartz space is def\/ined by
\begin{gather*}
\cS_{0,c}\big(\mR^{m}\big)=\bigcap_{k=0}^{\infty}D\big(L^{k}\big)
\end{gather*}
and one can check that the Fourier transform $\cF_{0,c}$ is an isomorphism of this space.
\end{remark}

\begin{remark}
In the limit case $c=0$, we can check that the kernel reduces to
\begin{gather*}
K(x,y)=\sum_{k=0}^{+\infty}\frac{k+\l}{\l}(-i)^{k}z^{-\l}J_{k+\l}(z)C_{k}^{\l}(w).
\end{gather*}
This is a~well-known expansion of the classical Fourier kernel (see~\cite[Section
11.5]{MR0010746}):
\begin{gather*}
K(x,y)=\frac{e^{-i\la\ux,\uy\ra}}{\Gamma(m/2)2^{\frac{m-2}{2}}}.
\end{gather*}
\end{remark}

We can now summarize the main properties of the deformed Fourier transform in the following theorem.

\begin{theorem}\label{StatementsFourier}
The operator $\cF_{0,c}$ defines a~unitary operator on $\cL^{2}_{0,c}(\mR^{m})$ and satisfies the
following intertwining relations on a~dense subset:
\begin{gather*}
\cF_{0, c} \circ \dD = i (1+c) \ux \circ \cF_{0, c},\qquad
\cF_{0, c} \circ \ux = \frac{i}{1+c}\dD \circ \cF_{0, c},\qquad
\cF_{0, c} \circ \mE = - \left( \mE + \delta \right) \circ \cF_{0, c}.
\end{gather*}
Moreover, $\cF_{0,c}$ is of finite order if and only if $c$ is rational.
\end{theorem}

\begin{proof}
Every $f$ in $\cL^{2}_{0,c}(\mR^{m})$ can be expanded in terms of the orthogonal basis
$\phi_{t,\ell,m}$, satisfying
\begin{gather*}
\langle\phi_{t_{1},\ell_{1},m_{1}},\phi_{t_{2},\ell_{2},m_{2}}\rangle=\delta_{t_{1}t_{2}}
\delta_{\ell_{1}\ell_{2}}\delta_{m_{1}m_{2}}\langle\phi_{t_{1},\ell_{1},m_{1}},
\phi_{t_{1},\ell_{1},m_{1}}\rangle,
\end{gather*}
see Section~\ref{section3}.
Note that the normalization can be computed explicitly (see~\cite[Theorem~6]{H12}).
As the eigenvalues of $\cF_{0,c}$ are given by (see~\eqref{eigvalsFourier})
\begin{gather*}
(-i)^{t}e^{-i\frac{\pi\ell}{2(1+c)}}
\end{gather*}
which clearly live on the unit circle, we conclude that
\begin{gather*}
\langle f,g\rangle=\langle\cF_{0,c}(f),\cF_{0,c}(g)\rangle
\end{gather*}
and that $\cF_{0,c}$ is a~unitary operator.

The intertwining relations are an immediate consequence of formula~\eqref{ActionOnBasis} combined
with the fact that $\phi_{t,\ell,m}$ is an eigenbasis of $\cF_{0,c}$.
The formula for $\mE$ follows from the anti-commutator (see Theorem~\ref{ospFamily})
\begin{gather*}
\{\dD,\ux\}=-2(1+c)\left(\mE+\frac{\delta}{2}\right).
\end{gather*}

The statement on the f\/inite order of the Fourier transform is an immediate consequence of the
explicit expression for the eigenvalues of the transform.
\end{proof}

Now we collect some properties of the kernel $K(x,y)$.

\begin{proposition}
One has, with $x,y\in\mR^{m}$
\begin{gather*}
K(\l x,y) = K(x, \l y), \quad \l > 0,\qquad
K(y,x) = \overline{K(x,y)},\\
K(0,y)  = \frac{1}{2^{\gamma_{0}/2-1}\Gamma(\gamma_{0}/2)},\qquad
K(\overline{s} x s,\overline{s} y s,) = \overline{s} K(x,y) s, \quad s \in \mathrm{Spin}(m),
\end{gather*}
where $\bar{\cdot}$ is the anti-involution on the Clifford algebra~$\cC l_{m}$.
\end{proposition}

\begin{proof}
The f\/irst property is trivial.
The second follows because
\begin{gather*}
\overline{\ux\wedge\uy}=-\ux\wedge\uy=\uy\wedge\ux.
\end{gather*}
The third property follows from Theorem~\ref{seriesrepr}.
Finally, the 4th equation follows because $z$ and $w$ are spin-invariant and
\begin{gather*}
\left(\overline{s}x s\right)\wedge\left(\overline{s}y s\right)=\overline{s}\left(\ux
\wedge\uy\right)s.
\tag*{\qed}
\end{gather*}
\renewcommand{\qed}{}
\end{proof}

We can also obtain Bochner identities for the deformed Fourier transform.
They are given in the following proposition.

\begin{proposition}\label{Bochner}
Let $M_{\ell}\in\cM_{\ell}$ be a~spherical monogenic of degree $\ell$.
Let $f(x)=f(r)$ be a~radial function.
Then the Fourier transform of $f(r)M_{\ell}$ and $f(r)\ux M_{\ell}$ can be computed as follows
\begin{gather*}
\cF_{0, c}(f(r) M_{\ell}) = e^{-\frac{i \pi \ell}{2(1+c)}} M_{{\ell}}(\uy') \int_{0}^{+\infty}
r^{\ell} f(r) z^{-\frac{\delta-2}{2}} J_{\frac{\gamma_{k}}{2} -1}(z) h(r) r^{m-1} dr,\\
\cF_{0, c} (f(r) \ux M_{\ell}) = -i e^{-\frac{i \pi \ell}{2(1+c)}} \uy' M_{{\ell}}(\uy')
\int_{0}^{+\infty} r^{\ell+1} f(r) z^{-\frac{\delta-2}{2}} J_{\frac{\gamma_{k}}{2} }(z)
h(r)r^{m-1} dr
\end{gather*}
with $\uy=s\uy'$, $\uy'\in\mS^{m-1}$ and $z=r s$.
\end{proposition}

\begin{proof}
This follows immediately from Theorem~\ref{seriesrepr} combined with Proposition \ref{reprkernels}.
\end{proof}

\begin{remark}
As a~special case of this proposition, we reobtain the eigenfunctions of the Fourier transform by
putting $f(r)=L_{t}^{\frac{\gamma_{\ell}}{2}-1}(r^{2})r^{\beta_{\ell}}e^{-r^{2}/2}$, resp.\
$f(r)=L_{t}^{\frac{\gamma_{\ell}}{2}}(r^{2})r^{\beta_{\ell}}e^{-r^{2}/2}$ (see equation~\eqref{eigvalsFourier}).
\end{remark}

Now we prove the following lemma.

\begin{lemma}\label{HeisenbergLemma}
For all $f\in\cL^{2}_{0,c}(\mR^{m})$ one has
\begin{gather*}
||\ux f(x)||^{2}+||\ux\left(\cF_{0,c}f\right)(x)||^{2}\geq\d||f(x)||^{2}.
\end{gather*}
The equality holds if and only if $f$ is a~multiple of~$e^{-r^{2}/2}$.
\end{lemma}

\begin{proof}
Using formula~\eqref{adjointsIP} and the unitarity of $\cF_{0,c}$, one can compute that
\begin{gather*}
||\ux f(x)||^{2}+||\ux\left(\cF_{0,c}f\right)(x)||^{2}=\frac{1}{(1+c)^{2}}\la
\left(\dD^{2}-(1+c)^{2}\ux^{2}\right)f,f\ra.
\end{gather*}
Now use the fact that the smallest eigenvalue of
\begin{gather*}
\frac{1}{(1+c)^{2}}\left(\dD^{2}-(1+c)^{2}\ux^{2}\right)
\end{gather*}
is given by $\d$, see Theorem~\ref{HarmOsc}.
This proves the inequality.

The equality holds when $f$ is a~multiple of an eigenfunction corresponding to the smallest
eigenvalue, i.e.\
when
$f$ is a~multiple of $e^{-r^{2}/2}$.
\end{proof}

This lemma allows us to obtain the Heisenberg inequality for the deformed Fourier transform
\begin{proposition}\label{Heisenberg}
For all $f\in\cL^{2}_{0,c}(\mR^{m})$, the deformed Fourier transform satisfies
\begin{gather*}
||\ux f(x)|| \cdot
||\ux\left(\cF_{0,c}f\right)(x)||\geq\frac{\d}{2}||f(x)||^{2}.
\end{gather*}
The equality holds if and only if~$f$ is of the form $f(x)=\l e^{-r^{2}/\a}$.
\end{proposition}

\begin{proof}
Using Lemma~\ref{HeisenbergLemma}, we can continue in the same way as in the proof of Theorem 5.28
in~\cite{Orsted2}.
\end{proof}

Now we can obtain the Master formula for the kernel of the Fourier transform.
We use the formula (see~\cite[p.~50]{MR0698780})
\begin{gather}
\label{besselint1}
\int_{0}^{+\infty}J_{\nu}(a t)J_{\nu}(b t)e^{-\gamma^{2}t^{2}}t dt=\frac{1}{2}\gamma^{-2}
e^{-\frac{a^{2}+b^{2}}{4\gamma^{2}}}I_{\nu}\left(\frac{a b}{2\gamma^{2}}\right),\qquad
\operatorname{Re}\nu>-1,\quad\operatorname{Re}\gamma^{2}>0,
\end{gather}
where $I_{\nu}(z)=e^{-i\frac{\pi\nu}{2}}J_{\nu}(i z)$.

We then obtain

\begin{theorem}[Master formula]\label{MasterFormula}
Let $s>0$.
Then one has
\begin{gather*}
\int_{\mR^{m}}K\left(y,x;i\frac{\pi}{2}\right)K\left(z,y;-i\frac{\pi}{2}\right)e^{-s r_{y}^{2}}h(r_{y})dy=\s_{m}
e^{-\frac{\o\d}{2}}K(z,x;\o)e^{-\frac{|\ux|^{2}+|\uz|^{2}}{2}\frac{1-\cosh{\o}}{\sinh{\o}}}
\end{gather*}
with $2s=\sinh{\o}$.
\end{theorem}

\begin{proof}
First observe that $K\left(y,x;i\frac{\pi}{2}\right)=K(y,x)$ and that $K(z,y;-i\frac{\pi}{2})$ is the
complex conjugate of $K\left(z,y;i\frac{\pi}{2}\right)$.

We rewrite the kernel $K$ obtained in Theorem~\ref{seriesrepr} in terms of the reproducing kernels~$P_{k}$ and~$Q_{k}$, i.e.\
as $K(x,y)=K_{0}(x,y)+K_{1}(x,y)$ with
\begin{gather*}
K_{0}(x,y) = \sum_{k=0}^{+\infty} \alpha_{k} \left( |\ux| |\uy|\right)^{-\frac{\delta-2}{2}}
J_{\frac{\gamma_{k}}{2} -1}(|\ux| |\uy|) P_{k}(\ux', \uy'),\\
K_{1}(x,y) = \sum_{k=0}^{+\infty} \beta_{k} \left( |\ux| |\uy|\right)^{-\frac{\delta-2}{2}}
J_{\frac{\gamma_{k}}{2}}(|\ux| |\uy|) Q_{k}(\ux', \uy'),
\end{gather*}
where $\alpha_{k}=e^{-\frac{i\pi k}{2(1+c)}}$ and $\beta_{k}=-i\alpha_{k}$.

When passing to spherical co-ordinates, the integral simplif\/ies, using Lemma~\ref{ReprKernOrth}, to
\begin{gather*}
\s_{m} \sum_{k=0}^{+\infty} \left( |\ux| |\uz|\right)^{-\frac{\delta-2}{2}} P_{k} (\uz', \ux')
\int_{0}^{+\infty} r e^{- s r^{2}} J_{\frac{\gamma_{k}}{2}-1} (r |\ux|)J_{\frac{\gamma_{k}}{2}-1}
(r |\uz|) dr\\
\qquad
{}+ \s_{m} \sum_{k=0}^{+\infty} \left( |\ux| |\uz|\right)^{-\frac{\delta-2}{2}} Q_{k} (\uz', \ux')
\int_{0}^{+\infty} r e^{- s r^{2}} J_{\frac{\gamma_{k}}{2}} (r |\ux|)J_{\frac{\gamma_{k}}{2}} (r
|\uz|) dr.
\end{gather*}
The radial integral can be computed explicitly using~\eqref{besselint1}.
Comparing with formula~\eqref{seriessemigroup} and Theorem~\ref{seriesholom} leads to the statement
of the theorem.
\end{proof}

\begin{remark}
For the Dunkl transform (see, e.g.,~\cite{MR1973996, MR2274972}) and for the Clif\/ford--Fourier
transform (see~\cite{DBXu}) one can compute even a~more general integral of the form
\begin{gather*}
\int_{\mR^{m}}K\left(y,x;i\frac{\pi}{2}\right)K\left(z,y;-i\frac{\pi}{2}\right)f(r_{y})h(r_{y})dy
\end{gather*}
with $f(r_{y})$ an arbitrary radial function of suitable decay.
This is done by using the addition formula for the Bessel function
\begin{gather*}
u^{-\l}J_{\l}(u)=2^{\l}\Gamma(\l)\sum_{k=0}^{\infty}(k+\l)(r^{2}|\ux||\uz|)^{-\l}J_{k+\l}(r
|\ux|)J_{k+\l}(r|\uz|)C^{\l}_{k}(\la\ux',\uz'\ra)
\end{gather*}
with $u=r\sqrt{|\ux|^{2}+|\uz|^{2}-2\la\ux,\uz\ra}$ instead of formula~\eqref{besselint1}.
Here, we cannot do that, as the orders of the Bessel functions do not match the order of the
Gegenbauer polynomials.
\end{remark}

\begin{remark}
Theorem~\ref{MasterFormula} is the starting point for the study of a~generalized heat equation, see,
e.g.,~\cite[Lemma 4.5(1)]{MR1620515} in the context of Dunkl operators.
\end{remark}

\section{Further results for the kernel}\label{section8}

In this section we will always be working in the non-Dunkl case, i.e.\ we put the multiplicity
function $\k=0$.
Theorem~\ref{seriesrepr} implies that the kernel of our deformed Fourier transform is a~function of
the type
\begin{gather*}
K(x,y)=f(z,w)+\ux\wedge\uy g(z,w)
\end{gather*}
with $f$, $g$ scalar functions of the variables $z=|\ux||\uy|$ and $w=\la\ux,\uy\ra/z$.
On the other hand, this kernel needs to satisfy the system of PDEs
\begin{gather*}
\dD_{y} K(x,y)= -i(1+c) K(x,y) \ux,\qquad
\left( K(x,y) \dD_{x} \right) = -i (1+c) \uy K(x,y),
\end{gather*}
as can be deduced from Theorem~\ref{StatementsFourier}.
In order to rewrite this system in terms of the variables~$z$,~$w$, we f\/irst observe that
\begin{gather*}
\upx f(z,w) = r^{-2} \ux z \partial_{z} f(z, w) + \left( z^{-1} \uy - r^{-2} \ux w
\right)\partial_{w}f(z, w),\\
\mE f(z, w) = z \partial_{z} f(z, w),\qquad
\upx \left( \ux \wedge \uy \right) = (1-m) \uy.
\end{gather*}
Using these identities, one obtains that the kernel is determined by the following 2 PDEs:
\begin{gather}
(m-1+c)g + (1+c) z \partial_{z}g +\frac{1}{z}\partial_{w} f + i(1+c)f -i(1+c) zw g=0,\nonumber\\
(1+c) z \partial_{z}f - w \partial_{w} f -c zw g -(1+c)z^{2} w \partial_{z} g + z( w^{2}-1)
\partial_{w}g+i(1+c) z^{2} g=0.\label{PDEkernel}
\end{gather}

\begin{remark}
Note that, contrary to the case of the classical Fourier transform and the Dunkl transform, where
the kernel is uniquely determined by the system of PDEs
\begin{gather*}
T_{j,x}K(x,y)=i y_{j}K(x,y),\qquad j=1,\ldots,m
\end{gather*}
this is not the case for the kernel of the radially deformed Fourier transform.
In fact, one can observe that there exist several dif\/ferent types of solutions of~\eqref{PDEkernel}.
This is discussed in detail in~\cite{DBNS} for a~similar system of PDEs in the context of the
so-called Clif\/ford--Fourier transform (see~\cite{DBXu}).
\end{remark}

Now we show that it is suf\/f\/icient to solve this system in dimension $m=2$ and $m=3$.
Recall that the kernel $K(x,y)$ is given in Theorem~\ref{seriesrepr}.
To know this kernel, it is hence suf\/f\/icient to know the series
\begin{alignat*}{3}
& A_{\lambda}=\sum_{k=0}^{+\infty}
\alpha_{k} (k + \l) J_{\frac{\gamma_{k}}{2}-1}(z)  C_{k}^{\l}(w), \qquad & &D_{\lambda}
 =\sum_{k=0}^{+\infty}   \alpha_{k-1} J_{\frac{\gamma_{k-1}}{2}}(z)  C_{k}^{\l}(w),& \\
& B_{\lambda}=\sum_{k=0}^{+\infty}
\alpha_{k} J_{\frac{\gamma_{k}}{2}-1}(z)  C_{k}^{\l}(w), \qquad & &E_{\lambda}
=\sum_{k=1}^{+\infty}   \alpha_{k}  J_{\frac{\gamma_{k}}{2}-1}(z)  C_{k-1}^{\l+1}(w),&\\
& C_{\lambda}=\sum_{k=0}^{+\infty}
\alpha_{k-1} (k + \l) J_{\frac{\gamma_{k-1}}{2}}(z)  C_{k}^{\l}(w), \qquad &&F_{\lambda}
=\sum_{k=1}^{+\infty}   \alpha_{k-1} J_{\frac{\gamma_{k-1}}{2}}(z)  C_{k-1}^{\l+1}(w),&
\end{alignat*}
because then one has
\begin{gather*}
K=\frac{1}{2\l}z^{-\frac{\d-2}{2}}\left(A_{\l}-i C_{\l}\right)+\frac{1}{2}
z^{-\frac{\d-2}{2}}\left(B_{\l}+i D_{\l}\right)-z^{-\frac{\d}{2}}\ux\wedge\uy\left(
E_{\l}+i F_{\l}\right).
\end{gather*}
Using the well-known property of the Gegenbauer polynomials $2\l C_{k-1}^{\l+1}(w)=\partial_{w}
C_{k}^{\l}(w)$, we observe the following recursion relations
\begin{alignat*}{5}
& A_{\l+1}= e^{i\frac{\pi}{2(1+c)}} \frac{1}{2 \l}\partial_{w} A_{\l}, \qquad &&
B_{\l+1}= e^{i\frac{\pi}{2(1+c)}} \frac{1}{2 \l}\partial_{w} B_{\l},\qquad &&
C_{\l+1}= e^{-i\frac{\pi}{2(1+c)}} \frac{1}{2 \l}\partial_{w} C_{\l},&\\
&D _{\l+1} = e^{-i\frac{\pi}{2(1+c)}} \frac{1}{2 \l}\partial_{w} D_{\l},\qquad &&
E _{\l} =\frac{1}{2 \l} \partial_{w} B_{\l},\qquad &&
F _{\l} =\frac{1}{2 \l} \partial_{w} D_{\l}. &
\end{alignat*}
We conclude that it suf\/f\/ices to know $A_{\l}$, $B_{\l}$, $C_{\l}$ and $D_{\l}$ for $\l=0,1/2$ or
$m=2,3$.
At this point, the problem of f\/inding explicit expressions for these functions for special values
of the deformation parameter $c$ is still open.

\appendix

\section{Properties of Laguerre and Gegenbauer polynomials}\label{appendixA}

The generalized Laguerre polynomials $L_k^{(\alpha)}$ for $k\in\mN$ are def\/ined as
\begin{gather*}
L_k^{(\alpha)}(t)=\sum_{j=0}^{k}\frac{\Gamma(k+\alpha+1)}{j!(k-j)!\Gamma(j+\alpha+1)}(-t)^j
\end{gather*}
and satisfy the orthogonality relation (when $\alpha>-1$)
\begin{gather*}
\int_{0}^\infty t^{\alpha}
L^{(\alpha)}_k(t)L^{(\alpha)}_l(t)\exp(-t)dt=\delta_{kl}\frac{\Gamma(k+\alpha+1)}{k!}.
\end{gather*}

The Gegenbauer polynomials $C^{(\alpha)}_k(t)$ are a~special case of the Jacobi polynomials.
For $k\in\mN$ and $\alpha>-1/2$ they are def\/ined as
\begin{gather*}
C_k^{(\alpha)}(t)=\sum_{j=0}^{\lfloor
k/2\rfloor}(-1)^j\frac{\Gamma(k-j+\alpha)}{\Gamma(\alpha)j!(k-2j)!}(2t)^{k-2j}
\end{gather*}
and satisfy the orthogonality relation
\begin{gather*}
\int_{-1}^1C_k^{(\alpha)}(t)C_l^{(\alpha)}(t)\big(1-t^2\big)^{\alpha-\frac{1}{2}}dt=\delta_{kl}\frac{\pi2^{1
-2\alpha}\Gamma(k+2\alpha)}{k!(k+\alpha)(\Gamma(\alpha))^2}.
\end{gather*}
One can prove that there exists a~constant $B(\a)$ such that
\begin{gather}\label{estimateGeg}
\sup_{-1\leq t\leq1}\left|\frac{1}{\a}C_k^{(\alpha)}(t)\right|\leq B(\a)k^{2\a-1},\qquad
\forall\, k\in\mN,
\end{gather}
see~\cite[Lemma 4.9]{Orsted2}.

The Bessel function $J_{\nu}(z)$ is def\/ined using the following Taylor series
\begin{gather*}
J_{\nu}(z)=\sum_{k=0}^\infty\frac{(-1)^k}{k!\Gamma(k+\nu+1)}
{\left({\frac{z}{2}}\right)}^{2k+\nu}.
\end{gather*}
For $z\in\mC$ and $\nu\geq-1/2$ one has the inequality (see, e.g.,~\cite{Sz})
\begin{gather}\label{estimateBessel}
\left|\left(\frac{z}{2}\right)^{-\nu}J_{\nu}(z)\right|\leq\frac{1}{\Gamma(\nu+1)}e^{|\operatorname{Im}{z}|}.
\end{gather}

\section{List of notations}\label{appendixB}

List of notations used in this paper:
\begin{alignat*}{3}
& m\quad &&\text{dimension of} \ \mR^{m},& \\
& \k \quad &&\text{multiplicity function on root system},& \\
& \mu\quad &&\text{Dunkl-dimension},& \\
& c\quad &&\text{deformation parameter of} \ \dD, & \\
& \o \quad &&\text{semigroup parameter with} \ \operatorname{Re}{\o}\geq0,& \\
& \upx \quad &&\text{ordinary Dirac operator},& \\
& \cD_{k} \quad &&\text{Dunkl Dirac operator},& \\
& \dD\quad &&\text{radially deformed Dirac operator},& \\
& \cF_{\dD}^{\o}\quad &&\text{exponential form of the holomorphic semigroup},& \\
& \cF_{0,c}^{\o}\quad &&\text{integral form of the holomorphic semigroup},& \\
& \cF_{\dD}\quad &&\text{exponential form of the Fourier transform},& \\
& \cF_{0,c}\quad &&\text{integral form of the Fourier transform}. &
\end{alignat*}

We also have the following def\/initions:
\begin{gather*}
\mu = m + 2 \sum_{\alpha \in R_+} \k_{\alpha},\qquad
\l = \frac{m-2}{2},\qquad
\sigma_{m} = 2 \pi^{m/2}/\Gamma(m/2),\qquad
\delta  = 1 + \frac{\mu-1}{1+c},\\
\beta_{\ell} = - \frac{c}{1+c}\ell, \quad \ell \in \mN,\qquad
\gamma_{\ell} = \frac{2}{1+c}\left( \ell + \frac{\mu-2}{2}\right) + \frac{c+2}{1+c}, \quad \ell
\in \mN.
\end{gather*}

\emph{Notations for variables.} Let $\ux$ and $\uy$ be vector variables in $\mR^{m}$.
Then we denote
\begin{gather*}
z = |\ux| |\uy|,\qquad
w  = \la \ux, \uy \ra /z.
\end{gather*}
When using spherical co-ordinates, we use $\ux=r\ux'$ with $\ux'\in\mS^{m-1}$, hereby
implicitly identifying a~vector in the Clif\/ford algebra with a~vector in $\mR^{m}$.

\subsection*{Acknowledgements}

H.~De Bie would like to thank E.~Opdam and J.~Stokman for valuable input and discussions during his
visit to the Korteweg--de Vries Institute for Mathematics in Amsterdam.
This visit was supported by a~Postdoctoral Fellowship of the Research Foundation~-- Flanders (FWO).
The last two authors would like to acknowledge support of the research grant GA \v{C}R P201/12/G028.

\pdfbookmark[1]{References}{ref}
\LastPageEnding
\end{document}